\documentclass[10pt,a4paper]{article}

\usepackage{graphics,epsfig,theorem,latexsym,amssymb,amsmath, amsfonts}
\usepackage{times,epic,eepic}

\usepackage{pifont}

\usepackage{cite}
\usepackage{url}

\usepackage{color}

\usepackage{newlfont}

\usepackage{yfonts}

\def\G1{\hbox{$\displaystyle{\mbox{\ding{172}}}$}}

\def\bd{\begin{description}}
\def\ed{\end{description}}

\def\beq{\begin{equation}}
\def\eeq{\end{equation}}
\def\bea{\begin{eqnarray}}
\def\eea{\end{eqnarray}}
\def\beas{\begin{eqnarray*}}
\def\eeas{\end{eqnarray*}}

\newtheorem{postulate}{Methodological Postulate}

\begin{document}

\title{Independence of the grossone-based infinity methodology from non-standard
analysis and comments upon logical fallacies in some texts asserting
the opposite }

\newcommand{\nms}{\normalsize}
\author{\\   Yaroslav D. Sergeyev\footnote{Yaroslav D. Sergeyev,
Ph.D., D.Sc., D.H.C., is Distinguished Professor at the University
of Calabria, Rende, Italy.
 He is also  Professor (a part-time contract) at the  Lobachevsky State University,
  Nizhni Novgorod, Russia and Affiliated Researcher at the Institute of High Performance
  Computing and Networking of the National Research Council of Italy, e-mail:
  {\tt
yaro@dimes.unical.it}}$\,\,^,$\footnote{The author thanks Prof.
Daniel Moskovich, Ben-Gurion University of the Negev, Beer-Sheva,
Israel for providing a preliminary list of logical fallacies present
in \cite{Flunks}.} }


\date{}

\maketitle

\begin{abstract}

This paper  considers non-standard analysis and   a recently
introduced computational me\-tho\-dology based on the notion of \G1
(this symbol is called \emph{gross\-one}). The latter approach was
developed with the intention to allow one to work with infinities
and infinitesimals numerically in a unique computational framework
and in all the situations requiring these notions.  Non-stan\-dard
analysis is a classical purely symbolic technique that works with
ultrafilters, external and internal sets, standard and non-standard
numbers, etc. In its turn, the \G1-based methodology does not use
any of these notions and proposes a more physical treatment of
mathematical objects separating the objects from tools used to study
them. It both offers a possibility to create new numerical methods
using infinities and infinitesimals in floating-point computations
and allows one to study certain mathematical objects dealing with
infinity more accurately than it is done traditionally. In these
notes, we explain that even though both methodologies deal with
infinities and infinitesimals, they are independent and represent
two different philosophies of Mathematics that are not in a
conflict. It is proved that texts
\cite{Flunks,Gutman_Kutateladze_2008,Kutateladze_2011} asserting
that the \G1-based methodology is a part of non-standard analysis
unfortunately  contain several logical fallacies. Their attempt to
show that the \G1-based methodology can be formalized within
non-standard analysis is similar to trying to show that
constructivism can be reduced to the traditional Mathematics.

 \end{abstract}

\section{Introduction}

The \G1-based infinity theory has been   introduced  in
\cite{Sergeyev,informatica,Lagrange} (see also recent surveys
\cite{EMS,UMI}) where this new numeral, \G1, called \emph{grossone},
and the related computational methodology have been described. This
new way of looking at infinity is not related either to Cantor's
cardinals and ordinals or to non-standard analysis of Robinson  or
to Levi-Civita field. As is well  known, there exists a variety of
mathematics: traditional, formalistic, intuitionistic, and other
philosophies of mathematics (see, e.g.,
\cite{Linnebo,Lolli_fil,Shapiro}). The \G1-based methodology
proposes one more way to look at Mathematics and Computer Science
following the example of Physics where tools used to observe objects
limit our possibilities of the observation.

 This methodological proposal has
attracted an appreciable amount of attention both of the pure and
applied mathematical communities. A number of papers studying
connections of the \G1-based approach to the historical panorama of
ideas dealing with infinities and infinitesimals (see
\cite{Lolli,Lolli_2,MM_bijection,Sorbi,Sergeyev_Garro}) have been
published. In particular, metamathematical investigations on the new
theory and its consistency can be found in \cite{Lolli_2}.  A number
of  reviews were published  in \textit{MIT Technology Review} and in
international scientific journals (see
\cite{Adamatzky,Pardalos,Prokopyev,Shylo,Trigiante}), etc. The
author received several international prizes and other
distinctions\footnote{Khwarismi International Award, assigned by The
Ministry of Science and Technology of Iran, 2017; Honorary
Fellowship, the highest distinction of the European Society of
Computational Methods in Sciences, Engineering and Technology, 2015;
Outstanding Achievement Award from the 2015 World Congress in
Computer Science, Computer Engineering, and Applied Computing, USA;
Degree of Honorary Doctor from Glushkov Institute of Cybernetics of
National Academy of Sciences of Ukraine, 2013; Pythagoras
International Prize in Mathematics, Italy, assigned by the city of
Crotone (where Pythagoras lived and founded his famous scientific
school) and the Calabria Region under the high patronage of the
President of the Italian Republic, Ministry of Cultural Assets and
Activities and the Ministry of Education, University and Research,
2010; Lagrange Lecture, Turin University, Italy, 2010; etc.} for
these results and was invited to present them as plenary lectures
and tutorials at more than 60 international congresses.

The \G1-based methodology   has   been successfully applied in
several areas of Mathematics and Computer Science: single and
multiple criteria optimization and ill-con\-di\-tio\-ning  (see
\cite{Cococcioni,DeLeone,DeLeone_2,Gaudioso&Giallombardo&Mukhametzhanov(2018),
homogeneity}), cellular automata (see \cite{DAlotto,DAlotto_3}),
Euclidean and hyperbolic geometry (see \cite{Margenstern_3}),
percolation (see \cite{Iudin_2,DeBartolo}), fractals (see
\cite{Caldarola_1,chaos,Menger,Biology,Koch,DeBartolo}),  infinite
series and the Riemann zeta function (see
\cite{Dif_Calculus,Riemann,UMI,EMS,Zhigljavsky}), the first Hilbert
problem and supertasks (see \cite{Rizza,first,EMS}), Turing machines
and probability (see
\cite{Rizza_2,EMS,Sergeyev_Garro,Sergeyev_Garro_2}), numerical
differentiation and   solution of  ordinary differential equations
(see \cite{ODE_3,Num_dif,ODE,ODE_2}), etc.

However, in the  paper \cite{Flunks} published in \emph{Foundations
of Science}\footnote{In Section 8, the authors of \cite{Flunks}
inform the reader  that before appearing in \emph{Foundations of
Science} their paper has been 5 times rejected by \emph{The
Mathematical Intelligencer.}} and in the papers
\cite{Gutman_Kutateladze_2008,Kutateladze_2011} published in two
journals printed by the Institute where   the authors of
\cite{Gutman_Kutateladze_2008,Kutateladze_2011} work, there are
numerous attacks on the \G1-based methodology and its author. The
paper \cite{Gutman_Kutateladze_2008} announces `a trivial
formalization of the theory of grossone' using non-standard
analysis. The paper \cite{Kutateladze_2011} does not have any
mathematical substance (not a single formula) and consists of
insults and accusations starting from the title. The paper
\cite{Flunks} in all seriousness attacks a playful note
\cite{medals} regarding several ways of counting Olympic medals won
by different countries and appeals   to
\cite{Gutman_Kutateladze_2008,Kutateladze_2011} trying again to show
that the \G1-based methodology can be reduced to non-standard
analysis.

Before going to   technicalities, let us compare   goals of the two
methodologies.  The successfully reached goal of the creator of
non-standard analysis, Abraham Robinson, was to reformulate
classical analysis and show that ideas of Leibniz can be put in a
form satisfying the requirements of rigor of the $XX^{th}$ century.
In fact, he wrote in paragraph 1.1 of his book \cite{Robinson}: `It
is shown in this book that Leibniz's ideas can be fully vindicated
and that they lead to a novel and fruitful approach to
\textit{classical analysis} and to many other branches of
Mathematics' (italics mine). In fact, classical analysis
reformulated using the language of Robinson was then used to address
such areas as topology, probability, etc. showing a notable
potential of this approach. However, as Davis writes in the
introduction to his book `Applied Nonstandard Analysis' (see
\cite{Davis_2}), `Nonstandard analysis is a technique rather than a
subject. Aside from theorems that tell us that nonstandard notions
are equivalent to corresponding standard notions, all the results we
obtain can be proved by standard methods.'

The goals of the \G1-based theory are different. The first goal is
to bridge the gap between modern Physics and Mathematics at least
partially, insofar as it concerns the separation of an object and a
tool used to observe it existing in Physics and often not present in
Mathematics. The difference between numbers and numerals is
emphasized and it is shown that numeral systems limit our
capabilities in describing numbers and other mathematical
objects\footnote{A \textit{numeral} is a symbol (or a group of
symbols) that represents a \textit{number} that is a concept. The
same number can be represented by different numerals. For example,
symbols `4', `four', `IIII', and `IV'  are different numerals, but
they all represent the same number.}. Traditionally,   various kinds
of infinity are perceived by people as \emph{distinct} mathematical
objects. The \G1-based methodology shows that this is not the case.
It argues that the objects -- infinite numbers -- are \emph{the
same}, they are just \emph{viewed in different ways} by different
theories and with   different accuracies. For example, it is shown
that in analysis and set theory, $\infty$ and~$\aleph_0$ are not two
different objects, they are two aggregative  images of infinite
quantities that people see through different `lenses' provided by
analysis and set theory.

The second goal of the \G1-based methodology is  to describe a
numeral system expressing infinities and infinitesimals that can be
used in \emph{all} occasions where we need to work with them. This
is done in order to have a situation similar to what we have with
finite numbers where numerals expressing them can be used in all the
occasions we need finite quantities. The viewpoint on infinities and
infinitesimals expressed in the \G1-based methodology does not
require a knowledge of cardinals, ordinals, ultrafilters, standard
and non-standard numbers, internal and external sets, etc. It avoids
a number of set theoretical paradoxes related to infinity and, in
general, introduces  a significant simplification in several fields
of mathematics related to infinities and infinitesimals.

Finally, the third goal consists of describing a   computational
methodology using the introduced numeral system and a  computational
device called the Infinity Computer (patented in USA and EU, see
\cite{Sergeyev_patent}) working \textit{numerically}\footnote{Recall
that numerical computations are performed with   floating-point
numbers that can be  stored in a computer memory. Since the memory
is limited, mantissa and exponent of these numbers can assume only
certain values and, therefore, the quantity and the form of numerals
that can be used to express floating-point numbers are fixed. Due to
this fact, approximations are required during computations with them
because an arithmetic operation with two floating-point numbers
usually produces a result that is not a floating-point number  and,
as a consequence, this result should be approximated by a
floating-point number. In their turn, symbolic computations are the
exact algebraic manipulations with mathematical expressions
containing \textit{variables that have not any given value}. These
manipulations are more computationally expensive than numerical
computations and only relatively simple codes can be elaborated in
this way. \label{foot numerical}} with finite, infinite, and
infinitesimal numbers in a unique computational framework and in
accordance with Euclid's Common Notion no.~5 `The whole is greater
than the part'. This gives the possibility to propose numerical
algorithms of a new type working \emph{in the same way} with numbers
that can have different infinite, finite, and infinitesimal parts.
Analogously, traditional computers do not consider as special cases
numbers having thousands and numbers hundreds of units. They
elaborate all of them in the same way.

Even from this brief introduction it can be seen that the goals and
the methodological platforms of the two approaches are very
different. In fact,  non-stan\-dard analysis is a purely symbolic
technique that works with ultrafilters, external and internal sets,
standard and non-standard numbers, etc. whereas the \G1-based theory
does not use any of these notions, is focused on numerical
computations, separates mathematical objects from tools used to
study them, and takes into account that, as in Physics, the tools
used limit our possibilities of the observation and determine the
accuracy of the observed results. Therefore, the  aim preset by the
authors of \cite{Flunks}, i.e., to prove that the \G1-based
methodology   is a part of non-standard analysis, is doomed to
failure. Such attempt is similar to trying to show that
constructivism can be reduced to  traditional Mathematics.

Let us now describe the structure of this paper. First,  several
logical fallacies\footnote{A logical fallacy (see
\cite{logical_fallacies}) is a flaw in reasoning. Logical fallacies
are like tricks or illusions of thought, and they are often very
sneakily used by politicians and the media to manipulate people.
This and the following footnotes explaining the meaning of fallacies
were taken from \cite{logical_fallacies}.} present in the papers
\cite{Flunks,Gutman_Kutateladze_2008,Kutateladze_2011} are listed
and analyzed in Section~2. The subsequent Section is dedicated to
some general considerations comparing the two methodologies under
scrutiny. Section~4 shows that weak numeral systems limit our
capabilities in measuring infinite sets. Traditional and \G1-based
numeral systems are compared. Then, since the authors of
\cite{Flunks} criticize not only the \G1-based methodology but also
a lexicographic rank described in \cite{medals} and propose their
own algorithm to count Olympic medals, the two methods are compared
briefly in Section~5. Finally, Section~6 concludes the paper.

\section{Logical fallacies committed  in   texts \cite{Flunks,Gutman_Kutateladze_2008,Kutateladze_2011}
}

In the following, five main logical fallacies present in
\cite{Flunks,Gutman_Kutateladze_2008,Kutateladze_2011} are
discussed.

\textbf{1. Confirmation Bias}. (\emph{We have a proclivity  to see
and agree with ideas that fit our preconceptions, and to ignore and
dismiss information that conflicts with them\footnote{You could say
that this is the mother of all biases, as it affects so much of our
thinking through motivated reasoning.   }}). The authors of
\cite{Flunks,Gutman_Kutateladze_2008,Kutateladze_2011} have an
expertise in non-standard analysis and try to defend it from the
\G1-based theory saying many times that `Sergeyev opposes his system
to nonstandard analysis' (see, e.g., \cite{Gutman_Kutateladze_2008},
page 1). In so doing   they seem not to notice that in his papers,
Sergeyev does not attack non-standard analysis. On the contrary, he
repeats again and again that his approach does not   oppose to any
existing methodology, it is just one more view on Mathematics.

The foundational platform of the \G1-based approach consists of
three methodological postulates and an axiom (the Infinite Unit
Axiom describing properties of grossone) that is added to axioms for
real numbers (see Appendix). In order to prove that Sergeyev's
methodology is a part of non-standard analysis it is necessary to
show that all the postulates and the axiom can be modeled using
non-standard analysis. However, Postulates~1 and~2 are not mentioned
at all in \cite{Flunks,Gutman_Kutateladze_2008}. The paper
\cite{Kutateladze_2011} lists the postulates but the only comment
related to them is the ironic exclamation `The scientific depth of
Sergeyev's postulates transpires'. There is no   other explanation
showing that the postulates can be (or were) considered by
non-standard analysis. The Infinite Unit Axiom is discussed in
\cite{Gutman_Kutateladze_2008}, but only partially.  The authors of
\cite{Gutman_Kutateladze_2008} write explicitly in page~2 that they
do not accept its main part   (this point will be debated below in
detail). The above said suffices to affirm that papers
\cite{Flunks,Gutman_Kutateladze_2008,Kutateladze_2011} have not
proved that Sergeyev's methodology is a part of non-standard
analysis. However, there are other fallacies to add.

{\bf 2. Personal Incredulity}. (\emph{Because you found something
difficult to understand, or are unaware of how it works, you made
out like it is probably not true\footnote{Complex subjects   require
some amount of understanding before one is able to make an informed
judgement about the subject at hand; this fallacy is usually used in
place of that understanding.}}). In Section~2 of \cite{Flunks}, its
authors write: `Sergeyev's attempted definition of \G1 as somehow
the number of elements of the set $\mathbb{N}$ contradicts other
passages where \G1 is included as a member of $\mathbb{N}$'. This
personal incredulity is a consequence of the previous fallacy. The
authors identify the set of natural numbers with its representation
offered by non-standard analysis. This representation does not
consider an entity being the number of elements of~$\mathbb{N}$.
Then, they apply non-standard analysis to the \G1-based theory
without noticing that methodological platforms of the two approaches
are different thus committing  a fallacy. Both approaches observe
the same set but do it using different mathematical tools (in
particular, different notations) and do not contradict one another.
The attempt to apply non-standard analysis tools in the \G1-based
framework is a fallacy similar to appealing to proofs by
contradiction in the frame of constructivism.

The Infinite Unit Axiom introduces \G1 as the number of elements of
the set $\mathbb{N}$ (see Appendix and \cite{EMS} for a detailed
discussion). This is performed by extrapolating from finite to
infinite the idea that $n$ is both the number of elements of the set
$\{ 1, 2, 3, \ldots , n-1, n \}$ and the last element of this set.
The fact that non-standard analysis does not offer this vision is
one more confirmation that the two theories are independent.

Another difficulty that the authors of \cite{Flunks} confess in
Section~4.1 in their appeal to works of Shamseddine (who, by the
way, twice delivered plenary lectures at conferences organized by
Sergeyev) is that they do not understand the difference between
numbers and numerals and between symbolic computations and numerical
ones.
 Even though Sergeyev in his papers dedicated
a lot of space to these topics, the above confession suggests that
it is necessary to return to them. These issues will be re-discussed
one more time below.

{\bf 3. The Fallacy Fallacy}. (\emph{Presuming a claim to be
necessarily wrong because a fallacy has been committed\footnote{It
is entirely possible to make a claim that is false yet argue with
logical coherency for that claim, just as it is possible to make a
claim that is true and justify it with various fallacies and poor
arguments.}}). Rough expository of Sergeyev's concepts in discussion
sections in \cite{Flunks,Gutman_Kutateladze_2008,Kutateladze_2011}
are seized upon and used to argue against grossone. At the same
time, the paper \cite{Flunks} does not contain a definition of
grossone anywhere in it, nor even its properties. It is never
explained what exactly is being attacked. A similar situation holds
in \cite{Gutman_Kutateladze_2008} where the Infinite Unit Axiom is
mentioned but methodological postulates are not provided. Vice
versa, \cite{Kutateladze_2011} writes about postulates but does not
mention the axiom. The impression is that the authors of
\cite{Flunks,Gutman_Kutateladze_2008,Kutateladze_2011}  avoid to
present the \G1-based theory in its complete form in order to make
their arguments more convincing.

{\bf 4. Appeal to Nature}. (\emph{Making the argument that because
something is `natural' it is therefore valid, justified, inevitable,
good, or ideal\footnote{Many `natural' things are also considered
'good', and this can bias our thinking; but naturalness itself does
not make something good or bad.}}). The authors write in
\cite{Gutman_Kutateladze_2008} that they do not accept that \G1 is
the number of elements of the set $\mathbb{N}$ `if for no other
reason than the fact that the set $\mathbb{N}$ of naturals (in the
popular sense of this fundamental notion) has no greatest element'.
Can `the popular sense' be the reason to refuse a mathematical
definition?

{\bf 5. Strawman}. (\emph{Misrepresenting someone's argument to make
it easier to attack\footnote{By exaggerating, misrepresenting, or
just completely fabricating someone's argument, it is much easier to
present your own position as being reasonable, but this kind of
deceitfulness  serves to undermine rational debate.}}). This is the
main fallacy. The paper \cite{Gutman_Kutateladze_2008} announces `a
trivial formalization of the theory of grossone' using non-standard
analysis. However, instead of working with \G1 introduced by
Sergeyev as the number of elements of the set $\mathbb{N}$, the
authors of \cite{Gutman_Kutateladze_2008} introduce a non-standard
object, declare that \emph{it is not} the number of elements of the
set $\mathbb{N}$, call the introduced object \emph{grossone}, use
the symbol \G1 to indicate it, and then proceed to rip it down. How
can this operation be called `formalization of the theory of
grossone'? Certainly, one can define a nonstandard object and study
it but how~is~this~related~to~what~Sergeyev~does?

\section{Comments emphasizing some other differences\\ between the \G1-based methodology and non-standard\\ analysis}

Let  us start this section by discussing  the object proposed to be
taken as grossone in \cite{Gutman_Kutateladze_2008} and, in addition
to the arguments provided above that are already sufficient to show
the independence of the two methodologies, explain why this concrete
proposal does not work. The authors of
\cite{Gutman_Kutateladze_2008} write at page 1:
\begin{quote}
`Fix an arbitrary infinitely large natural $\nu$ and denote its
factorial by \G1:
\end{quote}
\vspace{-5mm}
 \beq
\G1 = \nu!, \,\,\, \mbox{where}\,\,\, \nu \in \mathbb{N},\,\, \nu
\approx \infty.\mbox{'}
  \label{fal1}
 \eeq

In order to avoid the confusion present in papers
\cite{Flunks,Gutman_Kutateladze_2008,Kutateladze_2011}, let us use
hereinafter the symbol~\G1 for the original grossone introduced by
Sergeyev and the symbol $n(\nu)$ for the object used by the authors
of \cite{Gutman_Kutateladze_2008}, i.e., instead of (\ref{fal1}) we
write
 \beq
n(\nu) = \nu!, \,\,\, \mbox{where}\,\,\, \nu \in \mathbb{N},\,\, \nu
\approx \infty.
 \label{fal2}
 \eeq
 The authors of \cite{Gutman_Kutateladze_2008} correctly affirm that
$n(\nu)$ shares with \G1 \emph{Infinity} and \emph{Identity}
properties of the Infinite Unit Axiom  (see Appendix) and $n(\nu)$
is divisible by all finite numbers as \G1. However, this is not
sufficient to prove that $n(\nu)=\G1$. In fact, there are several
problems related to the introduction of $n(\nu)$ in (\ref{fal2}) as
a candidate for~\G1.

First, since $\nu$ is \emph{an arbitrary} infinitely large
non-standard natural number, $n(\nu)$ is not the number of elements
of $\mathbb{N}$ (in fact, at page~2 the authors of
\cite{Gutman_Kutateladze_2008} say this explicitly). However, it can
be seen from \emph{Divisibility} part of the Infinite Unit Axiom
 that, by taking $k=n=1$ in (\ref{3.3}) (see Appendix) we obtain that
the number of elements of $\mathbb{N}$ is equal to \G1. Since
$n(\nu)$ does not satisfy this condition, it cannot be chosen
as~\G1.

Second, the expression $\nu \approx \infty$ is meaningless in the
\G1-based framework since the traditional symbols $\infty, \aleph_0,
\omega, \aleph_1, $ etc. are not defined in it and, as a result,
cannot be used.

The third difficulty regards the words `Fix an arbitrary infinitely
large natural $\nu$'. The authors of \cite{Gutman_Kutateladze_2008}
do not explain how they intend to execute this operation of fixing.
In the \G1-based methodology,   \emph{fixing a variable} means that
a value is assigned to it, as it happens with variables in the
finite case. For instance, if we consider a finite $n$ then   we can
use, e.g., the numeral 34 and assign this value to it, i.e., $n=34$.
Analogously, for an infinite $n$, \G1-based numerals can be used and
we can assign values to $n$ using these numerals, e.g., $n=\G1-1$ or
$n=3\G1^2$, fixing so $n$.

In contrast, if we consider a non-standard infinite $\nu$ then it is
not clear which numerals consisting of a finite number of
symbols\footnote{Notice that the finiteness of the number of symbols
in the numeral is necessary for executing practical computations
since we should be able to write down and  store values we execute
operations with.} can be used to assign   a concrete value to $\nu$
since non-standard analysis does not provide numeral systems that
can be used for this purpose. In fact, all computations in
non-standard analysis theories are executed in a symbolic way
(recall footnote~\ref{foot numerical}) using a generic infinite
non-standard variable, e.g.,~$\nu$, and there is no possibility to
assign a value to $\nu$. Among other things, this means that
non-standard analysis is able to describe only those properties of
infinite numbers that are shared by all of them, since it does not
provide any instrument that would allow one to individualize an
infinite number, to distinguish it from another infinite number, and
to compare it with other infinite numbers. For instance, if one
considers two infinite non-standard numbers $\nu$ and $\xi$, where
$\nu$ is not expressed in terms of~$\xi$ (and vice versa), then it
is not clear how to compare them because of absence in non-standard
analysis of numeral systems that can express different values of
non-standard infinities. Notice that when we work with finite
quantities, then we can compare some $n$ and $k$ if they assume
numerical values, e.g., $k=25$ and $n=78$. Then, by using rules of
the numeral system the symbols 25 and 78 belong to, we can compute
that $n>k$. The same possibility is provided in the \G1-based
framework. For example, if $k=\G1-1$ and $n=3\G1^2$ then we can
compute
 \[
 n-k=3\G1^2-(\G1-1)= \G1(3\G1-1)+1  > 0
 \]
 and conclude that $n>k$.

It should be noticed here that this kind of  difficulties present in
non-standard analysis exists also in approaches dealing with
Levi-Chivita field,  since they work with \emph{a generic}
infinitesimal $\varepsilon$. Again it is not clear which numerals
can be used to assign a value to~$\varepsilon$ and to write
$\varepsilon=...$ (for instance, the nice web-based calculator
\emph{Inf} mentioned in \cite{Kutateladze_2011} and offered in
\cite{Crowell_Khafateh} operates with a generic symbol $d$ and there
is no possibility to assign a value to $d$). Moreover, approaches of
this kind leave unclear such issues as, e.g., whether the infinite
$1/\varepsilon$ is integer or not or whether $1/\varepsilon$ is the
number of elements of a concrete infinite set. The absence of
numeral systems allowing one to express quantities in non-standard
analysis and Levi-Civita field leads to a symbolic character of both
theories.

One more difference between $n(\nu)$ and \G1  can be indicated with
respect to infinite sequences.  An \textit{infinite sequence}
$\{a_n\}, a_n \in A, n \in \mathbb{N},$ is a function having as the
domain the set of natural numbers, $\mathbb{N}$, and as the codomain
a set $A$. A \textit{subsequence} is obtained from  a sequence by
deleting  some (or possibly none) of its elements.  In a sequence
$a_1, a_2, \ldots, a_{n-1}, a_n$ the number $n$ is the number of
elements of the sequence. Traditionally, only finite values of $n$
are considered. Grossone-based numerals give us the possibility to
observe infinite numbers and, therefore, to see not only the initial
elements of an infinite sequence $a_1, a_2, \ldots $ but also its
final part  $ \ldots, a_{n-1}, a_n$ where $n$ can assume different
infinite values. In other words, \G1-based theory allows one to
distinguish infinite sequences of different lengths. For instance,
the following two infinite sequences
\[
 \{a_n\} =  \{ 5,\hspace{3mm} 10,\hspace{3mm} \ldots \hspace{3mm} 5(\G1-1),\hspace{3mm}
5\G1 \},
\]
\[
 \{b_n\}  = \{  5,\hspace{3mm}10,\hspace{3mm} \ldots \hspace{3mm}
5 (\frac{2\G1}{5}-1),\hspace{3mm} 5\cdot \frac{2\G1}{5} \},
 \]
 have the same general element    $a_n=b_n=5n$ but they are
different because  the first sequence has \ding{172} elements and
the second sequence  has $\frac{2\G1}{5}$ elements.

Notice that since the set of natural numbers, $\mathbb{N}$,
has~\ding{172} elements, any sequence cannot have more than
\ding{172} elements (see \cite{EMS} for a detailed discussion). This
fact is very important in several research areas, in particular, in
that of  Turing machines (see \cite{Sergeyev_Garro}) where it allows
one to distinguish infinite tapes of different lengths. The
possibility to establish the maximal possible number of elements in
a sequence is not provided either by $n(\nu)$ or  by non-standard
analysis, in general.

\section{Numbers, numerals and measuring infinite sets }

Another important difference in comparison to non-standard analysis
consists of the fact that the \G1-based methodology through its
postulates 1 and 2 emphasizes the importance of numeral systems
being among our tools used to observe mathematical objects. Due to
postulate 1, the number of symbols we can use to write down numbers
is finite. Therefore, the choice of symbols and their meaning limit
our capabilities of observation of mathematical objects and
influence theoretical considerations, as well. The separation of
numbers (concepts) and numerals (symbols used to represent numbers)
was not discussed in depth traditionally and is not discussed in
non-standard analysis at all.

As an example, let us recall the Roman numeral system. It is not
able to express zero and negative numbers and such expressions as
III -- VIII or X -- X are indeterminate forms in this numeral
system. As a result, before the appearance of positional   systems
and the invention of zero (the second event was several hundred
years later with respect to the first one) mathematicians were not
able to create theorems involving zero and negative numbers and to
execute computations with them. The appearance of the positional
numeral system not only has allowed people to execute new operations
but has led to new theoretical results, as well. Thus, numeral
systems not only limit us in practical computations, they induce
constraints on reasoning in  theoretical considerations, as well.

Even a more significant, in the context of infinity, example of
limitations induced by numeral systems is provided by the numeral
system of Pirah\~{a}, a tribe living in Amazonia nowadays and
described in \textit{Science} in 2004 (see \cite{Gordon}). These
people use  an extremely simple numeral system for counting: one,
two, many. For Pirah\~{a}, all quantities larger than two are just
`many' and such operations as 2+2 and 2+1 give the same result,
i.e., `many'. Using their weak numeral system Pirah\~{a} are not
able to see, for instance, numbers 3, 4, and 5, to execute
arithmetical operations with them, and, in general, to say anything
about these numbers because in their language there are neither
words nor concepts for that\footnote{It should be noticed that the
astonishing numeral system of Pirah\~{a} is not an isolated example
of this way of counting. In fact, the same counting system, one,
two, many, is used by the Warlpiri people, aborigines living in  the
Northern Territory of Australia (see \cite{Butterworth}).}.

The poverty of the  numeral system of Pirah\~{a} leads   also to the
following results
 \beq \mbox{`many'}+ 1= \mbox{`many'},
\hspace{2mm}    \mbox{`many'} + 2 = \mbox{`many'},\hspace{2mm}
\mbox{`many'}+ \mbox{`many'} = \mbox{`many'} \label{1}
 \eeq
that are crucial for changing   our outlook on infinity. In fact, by
changing in these relations `many' with $\infty$ we get relations
used to work with infinity   in the traditional calculus
 \beq
\infty + 1= \infty,    \hspace{5mm}    \infty + 2 = \infty,
\hspace{5mm}   \infty +
  \infty =   \infty.
  \label{2}
 \eeq

Analogously, if we consider Cantor's cardinals (where, as usual,
numeral $\aleph_0$ is used for cardinality of denumerable sets and
numeral
 \textfrak{c}  for cardinality of the continuum) we
have similar relations
 \beq
 \aleph_0+  1 = \aleph_0, \hspace{5mm}    \aleph_0 + 2 =\aleph_0,
 \hspace{5mm}  \aleph_0+ \aleph_0 = \aleph_0,
  \label{3}
 \eeq
 \beq
 \textfrak{c} +  1 = \textfrak{c} , \hspace{5mm}    \textfrak{c}  + 2 = \textfrak{c} ,
 \hspace{5mm}  \textfrak{c}  + \textfrak{c}  = \textfrak{c} .
 \label{4}
 \eeq
 This comparison suggests that our difficulty in working with infinity is not a
consequence of the \textit{nature} of infinity but is a result of
\emph{weak numeral systems} having too little numerals  to express
the multitude of infinite numbers. As the possibility to have
numerals allowing us to express numbers 4, 5, etc. gives the
opportunity to execute more precise computations with finite numbers
in comparison with (\ref{1}), the introduction of the variety of
\G1-based numerals gives the possibility to execute more precise
w.r.t. (\ref{2})--(\ref{4}) computations with infinities (and
infinitesimals). Notice that this is not a new situation when the
introduction of   new symbols and   respective concepts allowed
mathematicians to have a progress in certain directions (it is
sufficient to mention $0, \infty, i, e,$ and $\pi$).

The separation of numbers from their representation is stressed in
the \G1-based theory whereas traditionally mathematicians do not pay
a particular attention to this issue. For instance, very often
mathematicians speak directly about real numbers and do not specify
which numeral systems are used to represent them. This can lead to
ambiguity in some cases.

 We illustrate
this statement by considering the following simple phrase `Let us
consider all $x \in [1,2]$'. For Pirah\~{a},
   \emph{all numbers} are just 1
 and 2. For people who do not know irrational numbers (or do not
 accept their existence) \emph{all numbers} are fractions $\frac{p}{q}$
 where $p$ and $q$ can be expressed in a numeral system they know.
If both $p$ and $q$ can assume values 1 and 2 (as it happens for
Pirah\~{a}), \emph{all numbers} in this case  are: 1,
$1+\frac{1}{2}$, and 2.
  For persons knowing positional numeral systems \emph{all numbers} are those
numbers that can be written in a positional system. Thus, in
different historical periods and in different cultures the phrase
`Let us consider all $x \in [1,2]$' has different meanings. As a
result, without fixing the numeral system we use to express numbers
we cannot fix the numbers we deal with and an ambiguity holds.

 In contrast, in
Physics, assertions are made with respect to what is visible at a
`lens' and not about the object located behind the lens. This is
done since  observations can be performed using different
instruments and without specifying the instruments assertions
regarding results have no meaning. For instance, the question: `What
do you see in this direction?'  is meaningless without indication a
tool used for the observation. In fact,  by eye the observer will
see certain things, by microscope other things, by telescope again
other things, etc.

\begin{table}[t]
 \caption{Cardinalities and the number of elements of some infinite sets.}
\begin{center} \small \label{table1}
\begin{tabular}{@{\extracolsep{\fill}}|c|c|c| }\hline
  Description     &  Cantor's  & Number  of \\
 of infinite sets  &   cardinalities  &   elements  \\
 \hline
  &   &     \vspace{-2mm}   \\
the set of natural numbers  $\mathbb{N}$   &  countable, $\aleph_0$ & \G1     \\
   &   &     \vspace{-2mm}   \\
  $\mathbb{N} \setminus \{ 3, 5, 10, 23, 114 \} $   &  countable, $\aleph_0$ & \G1-5     \\
  &   &     \vspace{-2mm}   \\
 the set of even numbers $\mathbb{E}$ (the set of odd numbers  $\mathbb{O}$)  &  countable, $\aleph_0$ & $\frac{\G1}{2}$    \\
 &   &     \vspace{-2mm}   \\
 the set of integers $\mathbb{Z}$    &  countable, $\aleph_0$ & 2\G1+1   \\
 &   &     \vspace{-2mm}\\
    $\mathbb{Z} \setminus \{ 0 \} $  &  countable, $\aleph_0$ & 2\G1    \\
 &   &     \vspace{-2mm}\\
  squares of natural numbers $\mathbb{G} = \{ x : x= n^2, x \in \mathbb{N},\,\, n   \in \mathbb{N} \}$   &  countable, $\aleph_0$ & $ \lfloor \sqrt{\G1} \rfloor$      \\
 &   &     \vspace{-2mm}   \\
   pairs of natural numbers $\mathbb{P}  = \{ (p,q) : p   \in
\mathbb{N},\,\, q \in \mathbb{N} \}$   &  countable, $\aleph_0$ &  $\G1^2$   \\
   &   &     \vspace{-2mm}\\
 the set of numerals  $\mathbb{Q}_1  = \{ \frac{p}{q}:   p   \in \mathbb{Z}, \,\, q \in \mathbb{Z},\,\,\, q \neq
 0
 \} $   &  countable, $\aleph_0$ &  $4\G1^2+2\G1$   \\
  &   &     \vspace{-2mm}\\
 the set of numerals  $\mathbb{Q}_2  = \{ 0,   -\frac{p}{q}, \,\, \frac{p}{q} : p   \in
\mathbb{N}, \,\,q \in \mathbb{N} \} $   &  countable, $\aleph_0$ &  $2\G1^2+1$   \\
  &   &     \vspace{-2mm}\\
  the power set of the set
  of natural numbers  $\mathbb{N}$   &  continuum, \textfrak{c}  & $2^{\mbox{\scriptsize{\ding{172}}}}$     \\
 &   &     \vspace{-2mm}   \\
 the power set of the set
  of even numbers  $\mathbb{E}$   &  continuum, \textfrak{c} & $2^{0.5\mbox{\scriptsize{\ding{172}}}}$     \\
 &   &     \vspace{-2mm}   \\
 the power set of the set
  of integers  $\mathbb{Z}$   &  continuum, \textfrak{c} & $2^{2\mbox{\scriptsize{\ding{172}}}+1}$     \\
 &   &     \vspace{-2mm}   \\
the power set of the set
   of numerals  $\mathbb{Q}_1$    &  continuum, \textfrak{c} & $2^{\mbox{\scriptsize{$4\G1^2+2\G1$}}}$     \\
 &   &     \vspace{-2mm}   \\
the power set of the set
   of numerals  $\mathbb{Q}_2$    &  continuum, \textfrak{c} & $2^{\mbox{\scriptsize{$2\G1^2+ 1$}}}$     \\
 &   &     \vspace{-2mm}   \\
  numbers $x \in [1,2)$ expressible
in the binary   numeral system  &  continuum, \textfrak{c} & $2^{\mbox{\scriptsize{\ding{172}}}}$     \\
 &   &     \vspace{-2mm}   \\
 numbers $x \in [1,2]$ expressible
in the binary   numeral system  &  continuum, \textfrak{c} & $2^{\mbox{\scriptsize{\ding{172}}}}+1$     \\
 &   &     \vspace{-2mm}   \\
numbers $x \in [1,2)$ expressible
in the decimal   numeral system   &  continuum, \textfrak{c} & $10^{\mbox{\scriptsize{\ding{172}}}} $    \\
 &   &     \vspace{-2mm}   \\
numbers $x \in [0,2)$ expressible
in the decimal   numeral system  &  continuum, \textfrak{c} & $2 \cdot 10^{\mbox{\scriptsize{\ding{172}}}}$    \\
\hline
\end{tabular}
\end{center}
\end{table}

The \G1-based theory follows physicists and does not talk, e.g.,
about the set of real numbers at the interval $[1,2)$ but about the
set of real numbers expressible in a fixed numeral system chosen to
represent real numbers over  $[1,2)$. Notice that Cantor's cardinals
do not allow us to distinguish the quantities  of real numbers over
$[1,2)$ written in the binary and in the decimal systems providing
the same answer: both sets have the cardinality of continuum. As
Table~\ref{table1} illustrates (see \cite{EMS} for its detailed
explanation), the \G1-based methodology allows us to register this
difference. In fact, the number of numerals expressing real numbers
over $[1,2)$ in the binary positional system is equal to
 $2^{\mbox{\scriptsize{\ding{172}}}} $ and the number of numerals expressing real numbers
over $[1,2)$ in the decimal positional system   is equal to
$10^{\mbox{\scriptsize{\ding{172}}}}$, where the number
$10^{\mbox{\scriptsize{\ding{172}}}}$  is infinitely larger than
$2^{\mbox{\scriptsize{\ding{172}}}}$. Moreover, sets of measure zero
are not present in the \G1-based framework and the accuracy of
measuring  infinite sets is equal to one element (i.e., if an
infinite set $A$ has $k$ elements where $k$ is expressed using
\G1-based numerals then exclusion/addition of one element from/to
$A$  gives the resulting set having exactly $k-1/k+1$ elements,
exactly as it happens with finite sets). This accuracy is
significantly higher than measuring executed by Cantor's cardinals.
For instance, for numerable sets one can see from Table~\ref{table1}
that excluding one number from $\mathbb{Z}$ can be registered by
\G1-based numerals. Analogously, it can be seen from the fourth and
the third lines from the end in Table~\ref{table1} that the number
of numerals expressed in the binary   numeral system and counted
over $[1,2)$ and $[1,2]$ are different.

To stress again the difference between numbers and numerals, notice
that the sets $\mathbb{Q}_1$ and $\mathbb{Q}_2$ in
Table~\ref{table1} are sets  of different rational numerals  and not
the sets of different rational numbers. For example, in the numeral
system $\mathbb{Q}_1$ the number 0 can be expressed by $2\G1$
different numerals
\[
\frac{0}{-\G1}, \,\, \frac{0}{-\mbox{\ding{172}+1}}, \,\,
\frac{0}{-\mbox{\ding{172}+2}}, \hspace{3mm} \ldots
\hspace{3mm}\frac{0}{-2}, \,\, \frac{0}{-1}, \,\, \frac{0}{1},  \,\,
\frac{0}{2},\hspace{3mm} \ldots
\hspace{3mm}\frac{0}{\mbox{\ding{172}-2}}, \,\,
\frac{0}{\mbox{\ding{172}-1}}, \,\, \frac{0}{\G1}.
\]
Analogously, e.g., such numerals as $\frac{-1}{-2}$, $\frac{1}{2}$,
and $\frac{2}{4}$ representing the same number have been calculated
as three different numerals.

The provided comparison between the traditional way to measure
infinite sets and the one provided by the \G1-based numerals shows
that the latter allow us to measure infinite sets more accurately
and gives one more argument to the impossibility of formalization of
the \G1-based methodology within non-standard analysis.

There are other  arguments to expose but what has been said is
already sufficient. The recent survey \cite{EMS} provides more
information and shows that \G1-based approach can be successfully
used in a variety of occasions where we need infinite and
infinitesimal quantities\footnote{For instance, \G1-based numerals
can be used for working with functions and their derivatives that
can assume different infinite, finite, and infinitesimal values and
can be defined over infinite and infinitesimal domains. The notions
of continuity and derivability can be introduced not only for
functions assuming finite values but for functions assuming infinite
and infinitesimal values, as well. Limits $\lim_{x \rightarrow a}
f(x)$ are substituted by expressions and $f(x)$ can be evaluated at
concrete   infinite or infinitesimal $x$ in the same way as it is
done with finite $x$. Series are substituted by sums having a
concrete infinite number of addends and for different number of
addends results (that can assume different infinite, finite or
infinitesimal values) are different as it happens for sums with a
finite number of summands. There are no divergent integrals, limits
of integration can be concrete different infinite, finite or
infinitesimal numbers and results  can assume different infinite,
finite or infinitesimal values. A number of set theoretical
paradoxes can be avoided, etc. }. On the one hand, \G1-based
numerals provide a higher accuracy of results with respect to
traditional symbols such as $\infty, \aleph_0, \omega$, etc. and, on
the other hand, eliminate the necessity to use different symbols in
different situations related to infinity. The fact that the same
symbols expressing infinities and infinitesimals can be used in all
the occasions we need them creates a unique framework with the
finite case, since the same numerals expressing finite quantities
are used in all the occasions we need them, as well.

In concluding this section, let us express a general doubt regarding
the meaning of attempts to reduce  the \G1-based methodology to
non-standard analysis. As was   discussed above, the \G1-based
methodology was developed to simplify the treatment of infinite sets
and to do away with any distinction between external and internal
sets, standard and non-standard numbers, cardinals and ordinals,
etc. In fact, it can be easily seen that the representation  of
infinities and infinitesimals in non-standard analysis using these
concepts is more tedious than the treatment   offered by the
\G1-based methodology (as was mentioned in \cite{Lolli}, the
\G1-based approach `... is simpler than non standard enlargements in
its conception, it does not require infinitistic constructions and
affords easier and stronger computation power.').

In their attempt to model \G1 using non-standard analysis the
authors of \cite{Flunks,Gutman_Kutateladze_2008,Kutateladze_2011}
not only commit fallacies but also make things significantly more
cumbersome with respect to Sergeyev's approach (e.g., in
\cite{Gutman_Kutateladze_2008} (see page 1 in bottom), its authors
try to render \G1 using external and internal sets departing
radically from Sergeyev's motivation). Traditionally, a complication
of a theory is acceptable if it produces more precise results.
Unfortunately, results provided in
\cite{Flunks,Gutman_Kutateladze_2008,Kutateladze_2011} are
\emph{less} precise. It is sufficient to mention that in
\cite{Gutman_Kutateladze_2008} (see the last page),  its authors
consider a set of numerals involving \G1 (where again both concepts,
numerals and grossone, do not coincide with their definitions used
by Sergeyev) and conclude that this set is countable. However, as
was already discussed above, \emph{countable} is a very rough result
with respect to the \G1-based approach since the introduction of \G1
allows one to measure infinite sets with the precision of one
element  (see Table~\ref{table1}).

Thus, what is the advantage to take a more precise and intuitive
theory and to try to reduce it to a less precise one that, in
addition, is more cumbersome?

\section{Once again on Olympic medals}
 The note \cite{medals}
entitled \textit{The Olympic Medals Ranks, Lexicographic Ordering,
and Numerical Infinities}   proposes a funny application of the
\G1-based theory related to counting Olympic medals. It is well
known that there exist several ways to rank countries with respect
to medals won during Olympics. If $g_A$ is the number of gold
medals, $s_A$ is the number of silver medals, and $b_A$ is the
number of the bronze ones won by   country~$A$, then many ranks have
the form
 \beq
R(\alpha,\beta,\gamma) = \alpha g_A + \beta s_A +  \gamma b_A,
  \label{fal3}
       \eeq
where $\alpha> 0, \beta> 0,$ and $\gamma> 0$ are certain weights.
However, the unofficial rank used by the Olympic Committee does not
allow one to use a numerical counter of the type (\ref{fal3}) for
ranking since it uses the lexicographic ordering to sort countries.
The rule applied for this purpose is the following: one gold medal
is more precious than any number of silver medals and one silver
medal is more precious than any number of bronze medals. The paper
\cite{medals} shows how it is possible to quantify  these words
\textit{more precious} by introducing a counter that for any a
priori unknown  finite number of medals allows one to compute a
numerical rank of a country using the number of gold, silver, and
bronze medals in such a way that the higher resulting number puts
the country in the higher position in the rank.   This can be easily
done by applying numerical computations with \G1-based numerals.

More formally, the problem considered in \cite{medals} is   to
introduce as a counter of  the type (\ref{fal3}) a number
$n(g_A,s_A, b_A)$ that should be calculated in such a way that for
countries $A$ and $B$ it follows that
 \beq
n(g_A,s_A, b_A) > n(g_B,s_B, b_B), \,\,  \mbox{if} \left\{
\begin{array}{l} g_A
>
g_B,\\
 g_A =  g_B,
s_A  > s_B,\\
g_A =  g_B, s_A  = s_B, b_A  > b_{B}.
\end{array} \right.
 \label{medals1}
       \eeq
In addition,   $n(g_A,s_A, b_A)$ should be introduced under
condition that the number $K > \max \{ g_A,s_A, b_A \}$ being an
upper bound for the number of medals of each type that can be won by
each country is unknown.

It is shown in \cite{medals} that  $n(g_A,s_A, b_A)$ can be
calculated using \G1 as follows
 \beq
n(g_A,s_A, b_A)=
 g_A \G1^{2} +   s_A \G1^1 +b_A \G1^0.   \label{medals4}
       \eeq
This formula gives us the rank of the type (\ref{fal3}) for the
country and this rank satisfies  condition (\ref{medals1}).  For
instance, let us consider the data
 \beq
  g_A=2, s_A=0, b_A=1, \hspace{1cm}
g_B=1,s_B=11, b_B=3.
 \label{medals3}
       \eeq
 Since \G1 is larger than any finite number (see the
Infinite Unit Axiom in Appendix), it follows from (\ref{medals4})
that
\[
n(g_A,s_A, b_A)= 2 \cdot \G1^{2} +  0 \cdot \G1^1 + 1 \cdot \G1^0 =
2 \G1^{2} +1
  >
\]
\[
n(g_B,s_B, b_B)= 1 \cdot \G1^{2} +  11 \cdot \G1^1 + 3 \cdot \G1^3 =
1   \G1^{2}  + 11   \G1^1 +3
\]
since
\[
2   \G1^{2} +1 - (1   \G1^{2}  +  11   \G1^1 + 3)= 1   \G1^{2} -11
\G1^1 -2= \underbrace{\G1 ( \G1 - 11)}_{\mbox{{\scriptsize
 positive and infinite}}}-2 > 0.
\]
 Thus,   to the set of
other existing ranks discussed in \cite{medals}   the grossone-based
rank (\ref{medals4}) has been added\footnote{Notice that the paper
\cite{medals} does not say that the rank (\ref{medals4}) is `the
best one'. It is just one more way to rank countries that can be
useful in certain situations.}. The paper \cite{medals} is concluded
by the obligatory comment that this way for counting ranks can   be
applied in all situations that require the lexicographic ordering
and not only for   three groups of objects but for any finite number
of them.

Let us now return to \cite{Flunks}. In section 7, its authors intend
to `demonstrate that the approach suggested by Sergeyev is useless'.
Unfortunately, in trying to do this, they once again show lack of
understanding of the meaning of the words `numerical computations'.
The authors of \cite{Flunks} write that Sergeyev regards
(\ref{medals4}) `as a ``numerical" rank just because it is a
``number" in the sense of his grossone theory'. This is not the
case. The   rank (\ref{medals4}) is numerical, since it allows one
to work on a computer (the Infinity Computer) where \G1-based
numbers are not manipulated as symbols and both $g_A, s_A, b_A$ and
the exponents 2, 1,   0 are expressed by floating-point numbers (see
footnote~\ref{foot numerical} to recall a brief explanation on the
essence of  numerical computations).

Then the  authors of \cite{Flunks} in Section 7 `indicate a very
simple and honest method' (are there dishonest methods?)  for
constructing a rank (\ref{fal3}) for lexicographic ordering. It
consists of using the binary numeral system  to compute the rank as
follows
 \beq
  r(g_A,s_A, b_A)= 0.
 \underbrace{11 \ldots\ldots 11}_{g_A \mbox{{\scriptsize
 positions}}}0\underbrace{11 \ldots\ldots 11}_{s_A
 \mbox{{\scriptsize
 positions}}}0\underbrace{11 \ldots\ldots 11}_{b_A \mbox{{\scriptsize
 positions}}}
.
 \label{fal4}
       \eeq
For the data (\ref{medals3}) this gives
\[
r(2,0,1) =0.11001 > r(1,11,3) =0.10111111111110111.
\]
However,   some comments can be done upon the very simple and honest
method proposed by the authors of \cite{Flunks}:
 \begin{enumerate}
 \item
 It does not assign the equal  weights to all medals of the same class.
Indeed, the authors of \cite{Flunks} explicitly state this because
in their opinion the first medal of a given class is more
significant achievement than the second one, that in its turn, is
more important than the third one, etc.\footnote{Notice that this
point of view  implies that the first competition is more important
than the second one, etc. violating so the principle of equality of
all sportive disciplines.} Thus, the very simple and honest method
does not address the same problem that is solved in \cite{medals}
(i.e., we face the Strawman fallacy again).
 \item
The very simple and honest method is not practical since   the
binary numeral system   requires  very long sequences of binary
digits to represent
 the rank (\ref{fal4}). In fact, even a very modest number of medals, e.g., 20
 gold, 20 silver, and 20 bronze, is already non representable in the
 IEEE 754 double-precision binary floating-point format having mantissa with 52 digits. In fact,
 the quantity of medals mentioned above requires 62  bits (60 positions for medals that
 are represented by 1 each and two zeros that separate groups of medals). Notice that the
 total
 number of medals in Olympics is usually significantly higher than
60. For instance, during The 2016 Summer Olympics  in Rio de Janeiro
there were 306 sets of medals. Since  countries can have several
athletes in each competition, and, therefore, a country can win
gold, silver, and bronze in each contest, 920 bits are required to
compute the rank proposed by the authors of \cite{Flunks}. This
means that a special data structure should be developed to implement
the very simple and honest method.
  \item
We recall that the whole story with the Olympics paper \cite{medals}
was conceived as  a divertissement.   However, the rank
(\ref{medals4}) proposed in \cite{medals} is not a joke,
  since it can be used in practical applications regarding the
  lexicographic ordering where \emph{the ordered quantities  are not necessary integer}. For
  example, in \cite{Cococcioni}, lexicographic multi-objective
linear programming problems have been considered, the gross-simplex
algorithm using the rank (\ref{medals4})  has been proposed and
implemented, results of numerical experiments have been provided.
The very simple and honest method proposed in \cite{Flunks} cannot
be used for this purpose since, due to (\ref{fal4}), it is able to
work with integer quantities $g_A, s_A,$ and $b_A$ only.
 \end{enumerate}

 Thus, the two ranks, (\ref{medals4}) and (\ref{fal4}), solve
 different problems, have different areas of applicability, and,
 therefore,
 are not in competition.

\section{A concluding remark}
The author of this paper hopes that this text will help people
interested in non-standard analysis to see that different points of
view on infinity are eligible.  The \G1-based approach does not
attack non-standard analysis in any way and there is no need to
defend it (especially, so inelegantly as it is done in
\cite{Flunks,Gutman_Kutateladze_2008,Kutateladze_2011}).
Non-standard analysis has been proposed in the middle of the
previous century, it has its vision of infinity, some scientists
practice it and obtain results in certain areas of Mathematics. This
is fine, no objection.

In its turn, the \G1-based methodology offers another vision of
infinity proposing a physically oriented Philosophy of Mathematics
and a number of related numerical algorithms. Some other scientists
use it in a variety of applications in pure and applied research
areas. In fact, nowadays there exist applications related to
numerical differentiation and optimization, ill-conditioning, ODEs,
traditional and blinking fractals, cellular automata, Euclidean and
hyperbolic geometry, percolation, probability, infinite series and
the Riemann zeta function, set theory and the first Hilbert problem,
Turing machines, etc.   (the interested reader is invited to consult
the recent survey \cite{EMS} and the dedicated web-page \cite{www}
for more information).

In conclusion, let us  quote a person who knew something about
infinity: `The essence of Mathematics lies entirely in its freedom.'
Certainly, the reader remembers that these words belong to Georg
Cantor.

\section*{Acknowledgement} The author thanks four unknown reviewers
for their valuable comments.

\section*{Appendix}

The \G1-based methodology is one of the possible views on infinite
and infinitesimal quantities and Mathematics, in general. It is
added to other existing philosophies of Mathematics such as
logicism, formalism, intuitionism, structuralism,  etc. (their
comprehensive analysis can be found, e.g., in
\cite{Linnebo,Lolli_fil,Shapiro}). Three postulates and an axiom
that is added to axioms for real numbers form the methodological
platform of the proposal. They are given below to make this paper
self-contained. A special attention in the \G1-based methodology is
paid to the fact than numeral systems that are among our tools used
to observe mathematical objects limit our capabilities of the
observation.   A detailed discussion on this methodological platform
can be found in \cite{EMS}.

\begin{postulate}
  \textit{We postulate  existence of infinite
and infinitesimal objects but accept that human beings and machines
are able to execute only a finite number of operations.}
\end{postulate}

\begin{postulate}
  \textit{We shall not   tell \textbf{what
are} the mathematical objects we deal with; we just shall construct
more powerful tools that will allow us to improve our capacities to
observe and to describe properties of mathematical objects.}
 \end{postulate}

\begin{postulate}
  \textit{We adopt the principle `The part is
less than the whole' to all numbers (finite, infinite, and
infinitesimal) and to all sets and processes (finite and infinite).}
\end{postulate}

\textbf{The  Infinite Unit Axiom.} The infinite unit of measure  is
introduced as the number of elements of  the set, $\mathbb{N}$, of
natural numbers. It is expressed by the numeral \ding{172} called
\textit{grossone} and has the following properties:

\textit{Infinity.}
  Any finite  natural number $n$  is less than grossone, i.e.,  $n
<~\G1$.

\textit{Identity.}
 The following
relations  link \ding{172} to identity elements 0 and 1
 \beq
 0 \cdot \G1 =
\G1 \cdot 0 = 0, \hspace{3mm} \G1-\G1= 0,\hspace{3mm}
\frac{\G1}{\G1}=1, \hspace{3mm} \G1^0=1, \hspace{3mm}
1^{\mbox{\tiny{\ding{172}}}}=1, \hspace{3mm}
0^{\mbox{\tiny{\ding{172}}}}=0.
 \label{3.2.1}
       \eeq

\textit{Divisibility.} For any finite natural number  $n$   sets
$\mathbb{N}_{k,n}, 1 \le k \le n,$ being the $n$th parts of the set,
$\mathbb{N}$, of natural numbers have the same number of elements
indicated by the numeral $\frac{\G1}{n}$ where
 \beq
 \mathbb{N}_{k,n} = \{k,
k+n, k+2n, k+3n, \ldots \}, \hspace{5mm} 1 \le k \le n, \hspace{5mm}
\bigcup_{k=1}^{n}\mathbb{N}_{k,n}=\mathbb{N}.
 \label{3.3}
       \eeq

\end{document}